\documentclass[
]{!!Arti}
\usepackage[leqno]{amsmath}
\usepackage[psamsfonts]{amssymb}
\usepackage{amsthm}
\usepackage[mathscr,psamsfonts]{eucal}
\usepackage{cmmib57}
\usepackage{!Style}
\usepackage{!Macros}
\usepackage{amscd}
\usepackage[nohug
]{diagrams}
\usepackage{color}
\usepackage{url}
\newarrow{Into}C--->
\newarrow{LR}{<}---{>}
\newcommand{\sten}{\,\ha\ten\,}
\newcommand{\rten}{\,\grave\otimes\,}
\newcommand{\lten}{\,\acute\otimes\,}
\DeclareMathOperator{\sdim}{{s-dim}}
\DeclareMathOperator{\Slin}{{{s-Lin}}}
\begin{document}
\title{\bf Semi--vector spaces\\ and units of measurement}
\author{
{\bf Josef Jany\v{s}ka$^1$, Marco Modugno$^2$, Raffaele Vitolo$^3$}
\bigskip
\\
\fz $^1$ Department of Mathematics and Statistics, Masaryk University
\\
\fz Jan\'a\v{c}kovo n\'am 2a, 602 00 Brno, Czech Republic
\\
\fz email: {\tt janyska@math.muni.cz}
\myskip
\\
\fz $^2$ Department of Applied Mathematics, Florence University
\\
\fz Via S. Marta 3, 50139 Florence, Italy
\\
\fz email: {\tt marco.modugno@unifi.it}
\myskip
\\
\fz $^3$ Department of Mathematics ``E. De Giorgi''
\\
\fz Via per Arnesano, 73100 Lecce, Italy
\\
\fz email: {\tt raffaele.vitolo@unile.it}
}
\date{}
\edition{\sm \emph{Preprint: 2007.07.26. - 17.20.}}
\pagestyle{headings}
\maketitle
\begin{abstract}
 This paper is aimed at introducing an algebraic model for physical
scales and units of measurement.
 This goal is achieved by means of the concept of ``positive space"
and its rational powers.
 Positive spaces are 1--dimensional ``semi--vector spaces'' without the
zero vector.
 A direct approach to this subject might be sufficient.
 On the other hand, a broader mathematical understanding requires the
notions of sesqui and semi--tensor products between semi--vector
spaces and vector spaces.

 So, the paper is devoted to an original contribution to the
algebraic theory of semi--vector spaces, to the algebraic analysis of
positive spaces and, eventually, to the algebraic model of physical
scales and units of measurement in terms of positive spaces.

\par\nopagebreak
\end{abstract}
\myintro
\label{Introduction}
 Most objects in physics (like metric field, electromagnetic field,
gauge fields, masses, charges, and so on) have physical dimensions,
as their numerical values depend on the choice of units of
measurement.
 Usually, these physical dimensions are dealt with in an ``informal''
way, without a precise mathematical setting.

 In recent years, a geometric formulation of covariant classical and
quantum mechanics includes also a formal mathematical setting of
physical scales and units of measurement in an original algebraic way
(see, for instance,
\cite{CanJadMod95,JanMod02c,JanMod05p2,ModSalTol05,SalVit00,Vit99,
Vit00}).
 Such a rigorous mathematical setting of scales plays an essential
role in some aspects of the covariant theory, in particular in the
classification of covariant operators (see, for instance,
\cite{JanMod02c}).

 This approach is based on the notion of ``positive space'' and
its rational powers, as model for the spaces of scales and units of
measurement.
 Then, tensor products between positive spaces and vector bundles
arising from spacetime yield ``scaled objects'', i.e. objects with
physical dimension.

 In the above papers, this mathematical scheme is sketched very
briefly.
 So, a comprehensive mathematical foundation of this subject is
required.
 
 ``Positive spaces'' are 1--dimensional ``semi--vector'' spaces.
 In principle, we could develop a formalism for positive spaces
directly, but a broader mathematical understanding of this subject
suggests to insert the theory of positive spaces in the wider
framework of semi--vector spaces.
 This concept is not new, and has been used in very different
contexts: in \cite{KKB72} for the analysis of some properties of
$\B Z_2$-valued matrices, in \cite{GahGah99} for problems of fuzzy
analysis, in \cite{Pap80} for problems of measure theory, in
\cite{PraSer74,PraSer76} for topological fixed point problems.
 Here, we introduce the sesqui and semi--tensor products
between semi--vector spaces and vector space and prove several
results as well.
 Such results are new and of independent interest in the algebraic
theory of semi--vector spaces.

 Thus, the goal of the paper is two folded: giving a contribution to
the algebraic theory of semi--vector spaces, with special emphasis to
tensor products, and proposing an algebraic model for physical
scales and units of measurement.

\myskip

 Thus, the paper is organised as follows.

 In section 1, we discuss semi--vector spaces over the semi--field
$\Rn^+ \,,$
introduce the sesqui and semi-tensors products between semi--vector
spaces and vector spaces.

 In section 2, we restrict our attention to positive spaces, which are
one-dimensional semi--vector spaces without zero vector, and introduce
their rational powers.

 In section 3, we discuss the algebraic model of physical scales and
units of measurement. 
 We start by assuming three positive spaces, 
$\B T \,,\, \B L$ 
and 
$\B M$ 
as representatives of the spaces of time, length and mass scales.
 Then, we describe all possible derived scales in terms of
semi--tensor products of rational powers of 
$\B T \,,\, \B L$
and
$\B M \,.$ 
 Next, we introduce the ``scaled objects'', by considering the
sesqui--tensor product of a positive space with a vector bundle
arising from spacetime.
 We also sketch differential calculus with scaled objects.
 Several of the above results are related to the dimensional analysis; we
briefly analyse the interplay between that theory and our algebraic setting.
\section{Semi--vector spaces}
\label{section: Semi--vector spaces}
\subsection{Semi--vector spaces}
\label{Semi--vector spaces}
\bsm
 We define the semi--vector spaces analogously to vector spaces, by
substituting the field of real numbers with the semi--field of
positive real numbers.
 Semi--vector spaces are similar to vector spaces in some respects;
however, the lack of some standard properties yields subtle
problems, which require a careful treatment.
\esm

 Let 
$\Rn^+ \sub \Rn$
and
$\Rn^- \sub \Rn$
be the subsets of positive and negative real numbers; moreover, let
us set
$\Rn^+_0 \byd \Rn^+ \cup \{0\}$
and
$\Rn^-_0 \byd \Rn^- \cup \{0\} \,.$

 The set
$\Rn^+$
is a \emp{semi--field} \cite{HebWei93} with respect to the
operations of addition and multiplication.
This means that 
$(\Rn^+,+)$ 
is a commutative semi--group,
$(\Rn^+,\cdot)$ 
is a commutative group and the distributive law holds.
 Note that, for each 
$x \,, y \,, z \in \Rn^+ \,,$ 
the ``cancellation law'' holds:
if
$x + z = y + z \,,$
then
$x = y \,.$

\bDf
 A \emp{semi--vector space} (over $\Rn^+$) is defined to be a set
$\f U$
equipped with the operations
$+ : \f U \car \f U \to \f U$
and
$\cdot : \Rn^+ \car \f U \to \f U \,,$
which fulfill the following properties, for each
$r,s \in \Rn^+, \, u,v,w \in \f U \,,$
\bgt
u + (v + w) = (u + v) + w \,,
\qquad
u + v = v + u \,,
\\
(r s) \, u = r \, (s u) \,,
\qquad
1 \, u = u \,,
\qquad
r \, (u + v) = r u + r v \,,
\qquad
(r + s) \, u = r u + s u \,.\END
\end{gather*}
\eDf

 Despite the fact that some authors (see \cite{GahGah99}) require that
a semi--vector space has a zero vector, here we do not make such a
general assumption.
 However, interesting properties arise for semi-vector spaces with a
zero vector.

\myskip

 Let
$\f U$
be a semi--vector space.

 An element
$0 \in \f U$
is said to be \emp{(additively) neutral} (or, equivalently, a
\emp{zero vector}) if, for each
$u \in \f U \,,$
we have
$0 + u = u \,.$

 If two elements
$0, 0' \in \f U$
are neutral, then they coincide; in fact, we have
$0' = 0 + 0' = 0' + 0 = 0 \,.$

 If
$0 \in \f U$
is a neutral element, then, for each
$r \in \Rn^+ \,,$
we have
$r 0 = 0 \,;$
in fact, for each
$u \in \f U \,,$
we have
$r 0 + u = r (0 + \fr1 r u) = r (\fr1 r u) = u \,.$

\myskip

 A semi--vector space equipped with a zero vector is said to be
\emp{complete}.

 A complete semi--vector space
$\f U$
turns out to be also a semi-vector space over
$\Rn^+_0 \,,$
by setting
$0 \, u = 0 \,,$
for each
$u \in \f U \,.$
	But we shall always refer to the scalar multiplication with respect
to the semi--field
$\Rn^+$
even for complete semi--vector spaces.

\myskip

 Let
$\f U$ 
be a complete semi--vector space.
 We say that 
$u \in \f U$ 
is \emp{invertible} if there exists a vector 
$v \in \f U$ 
such that 
$u+v = 0 \,.$

 A non complete semi--vector space, or a complete semi--vector space
with no invertible elements, is said to be \emp{simple}.

\myskip

 We say that a semi--vector space
$\f U$
is \emp{regular} if the
``cancellation law'' holds, i.e. if, for each
$u, u', v \in \f U \,,$ 
the equality
$u + v = u' + v$
implies
$u = u' \,.$

 In regular complete semi--vector spaces further properties hold.

\bNt
 Let
$\f U$
be regular and complete.
 If
$0' \in \f U$
is an element such that
$0' + u = u$
for a certain
$u \in \f U \,,$
then
$0' = 0 \,;$
in fact, we have
$0' + u = u = 0 + u \,,$
hence, the cancellation law yields
$0' = 0 \,.$
 Moreover, if
$u \in \f U$
is invertible, then its inverse is unique; in fact,
$u + v = 0 = u + v'$
implies
$v = v' \,,$
in virtue of the cancellation law.\END
\eNt

\bNt
 Let
$\f V$
be a vector space and
$\f U \sub \f V$
a subset, which turns out to be a semi--vector space with respect to
the sum and product inherited from
$\f V \,.$
 Then,
$\f U$
is regular.\END
\eNt

\bDf
 Let us consider a simple semi--vector space
$\f U \,.$

 A subset
$\f B \sub \f U$
is said to be a \emp{semi--basis} if any non neutral
$u \in \f U$
can be written in a unique way as
$u = \sum_{i \in I_u} u^i \, b_i \,,$
with
$u^i \in \Rn^+$
and
$b_i \in \f B \,,$
where the sum is extended to a finite family of indices
$I_u \,,$
which is uniquely determined by
$u \,.$

 Accordingly, we denote by
$\f B_u$
the finite subset
$\f B_u \byd \{b_i\}_{i \in I_u} \sub \f B \,,$
which is uniquely determined by
$u \,.$

 The semi--vector space
$\f U$
is said to be \emp{semi--free} if it admits a semi--basis.\END
\eDf

 We can easily see that if
$\f U$
is a complete semi--free semi--vector space and 
$\f B$
is a semi--basis, then
$0 \nin \f B \,.$

\bPr
 If
$\f U$
is semi--free, then it is regular.
\ePr

\bpf
 Let
$\f B$
be a semi--basis.
 If
$u,v,w \in \f U \,,$
then we have unique expressions of the type
$u = \sum_i u^i \, b_i \,,\; 
 v = \sum_j v^j \, b_j \,,\;
 w = \sum_h w^h \, b_h \,,$
with
$u^i, v^i, w^i \in \Rn^+$
and
$b_i \in \f B_u \,,
 b_j \in \f B_v \,,
 b_h \in \f B_w \,.$

 We have 
$\f B_{u+w} = \f B_u \cup \f B_w$ 
and
$u+w = \sum_k (\hat u^k + \hat w^k)b_k \,,$
where 
$b_k \in \f B_{u+w} \,,$ \,
$\ha u^k = u^k \,,$ 
or 
$\ha w^k = w^k \,,$ 
if 
$b_k \in \f B_u \,,$ 
or 
$b_k \in \f B_w \,,$
respectively, and 
$\hat u_k = 0 \,,$ 
or 
$\ha w_k = 0 \,,$ 
if 
$b_k \nin \f B_u \,,$ 
or 
$b_k\not\in\f B_w \,,$
respectively.
 Analogous considerations hold for
$v + w \,.$

Then, from the uniqueness of the expression of semi--vectors, if
$u + w = v + w \,,$
then we have 
$\f B_{u+w} = \f B_{v+w}$ 
and 
$\sum_k (\ha u^k + \ha w^k) \, b_k = 
 \sum_k (\ha v^k + \ha w^k) \, b_k \,,$
which implies 
$\ha u^k = \ha v^k \,,$
i.e. 
$u^k = v^k \,.$
 Hence,
$u = v \,.$\QED
\epf

\bPr\label{Proposition: bijection of bases}
 Let
$\f U$
be a semi--free semi--vector space.
 Moreover, let
$\f B$
be a semi--basis and
$\baf B \sub \f U$
a subset.
 Then, the following facts hold.

1) Let
$\baf B$
consist of non vanishing elements and suppose that there is a
bijection
$\baf B \to \f B : \ba b_i \mto b_i$
and, for each
$\ba b_i \in \baf B \,,$
we have 
$\f B_{\ba b_i} = \{b_i\} \,,$
so that we obtain the expression 
$\ba b_i = \ba S_i \, b_i \,,$
with
$\ba S_i \in \Rn^+ \,.$
 Then, 
$\baf B$
is a semi--basis and, for each
$b_i \in \f B \,,$
we have 
$\baf B_{b_i} = \{\ba b_i\} \,,$
so that we obtain the expression 
$b_i = S_i \, \ba b_i = (1/\ba S_i) \, \ba b_i \,.$

2) Suppose that
$\baf B$
be a semi--basis.
 Then, there is a bijection
$\baf B \to \f B : \ba b_i \mto b_i$
and, for each
$\ba b_i \in \baf B \,,$
we have 
$\f B_{\ba b_i} = \{b_i\} \,,$
so that we obtain the expression 
$\ba b_i = \ba S_i \, b_i \,,$
with
$\ba S_i \in \Rn^+ \,;$
moreover, for each
$b_i \in \f B \,,$
we have 
$\baf B_{b_i} = \{\ba b_i\} \,,$
so that we obtain the expression 
$b_i = S_i \, \ba b_i = (1/\ba S_i) \, \ba b_i \,.$
\ePr

\bpf
1) For each
$b_i \in \f B \,,$
the expression 
$\ba b_i = \ba S_i \, b_i \,,$
implies
$b_i = (1/\ba S_i) \, \ba b_i \,.$
 Then, for each
$0 \neq v \in \f U \,,$
we have
$v = \sum_{i \in I_v} v^i \, b_i = 
\sum_{j \in I_v} v^j \, (1/\ba S_j) \, \ba b_j \,.$
 Indeed, the 2nd expression is unique.
In fact, if
$v = \sum_{h \in \ba I_v} \ba v^h \, \ba b_h \,,$
then we have also
$v = \sum_{h \in \ba I_v} \ba v^h \, \ba S_h \, b_h \,;$
hence, in virtue of the uniqueness of the expression with respect to
the semi--basis
$\f B \,,$
we obtain
$I_v = \ba I_v$
and
$\ba v^i \, \ba S_i = v^i \,,$
i.e.
$\ba v^i = v^i \, (1/\ba S_i) \,.$

2) For each
$b_i \in \f B$
and
$\ba b_j \in \baf B \,,$
we have unique expressions of the type
\beq
b_i = \sum_{h \in \ba I_{b_i}} S^h_i \, \ba b_h
\ssep{and}
\ba b_j = \sum_{k \in I_{\ba b_j}} \ba S^k_j \, b_k \,,
\ssep{with}
S^h_i \,,\, \ba S^k_j \in \Rn^+ \,,
\eeq

 Then, by substituting the 2nd expressions into the 1st one, we
obtain 
\beq
b_i = \sum_{h \in \ba I_{b_i}} \sum_{k \in I_{\ba b_h}}
S^h_i \, \ba S^k_h \, b_k \,.
\eeq

 On the other hand, by recalling the uniqueness of the decomposition
of each element
$b_i \in \f B$
with respect to this semi--basis
$(b_i)$
and observing that all coefficients
$S^j_i$
and
$\ba S^h_j$
are positive, we deduce that the above equality holds if and only if
the following conditions hold:

a) there is a bijection
$\f B \to \baf B : b_i \mto \ba b_i \,,\;$

b)
$\baf B_{b_i} = \{\ba b_i\} \,,\,$
$\f B_{\ba b_j} = \{b_j\} \,,\;$

c)
$\ba b_i = \ba S_i \, b_i$
and
$b_i = S_i \, \ba b_i = (1/\ba S_i) \, \ba b_i \,,$
with
$\ba S_i \in \Rn^+ \,.$\QED
\epf

 We stress that the transition law between semi--bases of semi--free
semi--vector spaces given by the above 
Proposition \ref{Proposition: bijection of bases}
is essentially more restrictive than the transition law between bases
of vector spaces. 

\bCr
 If
$\f U$
is semi--free, then all semi--bases have the same cardinality.\END
\eCr

\bDf
 If
$\f U$
is a semi--free semi--vector space, then we define the
\emp{semi--dimension} of
$\f U$
to be the cardinality of its semi--bases. 
 We denote the semi--dimension of
$\f U$
by
$\sdim \f U \,.$\END
\eDf

 In the following, we shall refer only to semi--free semi-vector
spaces with finite semi--dimension in those formulas where we compute
explicitly the semi--dimension.
\subsection{Examples of semi--vector spaces}
\label{Examplers of semi--vector spaces}
\bsm
 We have examples of quite standard semi--vector spaces and odd
examples as well.
\esm

\bEx
 A vector space (over
$\Rn$) turns out to be a complete and regular semi--vector space.
 But, it is not simple and semi--free.\END
\eEx

 A semi--vector space is said to be \emph{proper} if it is not a
vector space.

\bEx
 If
$\f V$
is a vector space, then the subset
$\f U \byd \f V - \{0\} \sub \f V$
is not a semi--vector space, because the sum of an element and of its
negative is not defined in
$\f U \,.$\END
\eEx

\bNt
 Let
$\f V$
be a vector space and
$\f S \sub \f V$
a subset consisting of
$n \geq 1$
vectors.

 Then, we define the \emp{cone generated by}
$\f S$
\emp{over}
$\Rn^+$
to be the subset
$\la \f S \ra_+ \sub \f V$
consisting of all linear combinations of the type
$r^1 s_1 + \dots + r^m s_m \,,$
with
$m = 1, \dots, n$
and
$s_1, \dots, s_m \in \f S \,,\,$
$r^1, \dots, r^m \in \Rn^+ \,.$

  Analogously, we define the \emp{cone generated by}
$\f S$
\emp{over}
$\Rn^-$
to be the subset
$\la \f S \ra_- \sub \f V$
consisting of all linear combinations of the type
$r^1 s_1 + \dots + r^m s_m \,,$
with
$m = 1, \dots, n$
and
$s_1, \dots, s_m \in \f S \,,\,$
$r^1, \dots, r^m \in \Rn^- \,.$

 Clearly, if
$\f S_- \sub \f V$
is the subset consisting of the the negatives of the set
$\f S \,,$
then we have
$\la \f S_- \ra_+ = \la \f S \ra_-$
and
$\la \f S_- \ra_- = \la \f S \ra_+ \,.$\END
\eNt

\bNt
 Let
$\f V$
be a vector space and
$\f S \sub \f V$
a subset consisting of
$n \geq 1$
vectors. 
 The set 
$\la \f S \ra_+$ 
is convex in 
$\f V \,.$
\eNt

\bpf
 The set  
$\la \f S \ra_+$ 
is convex if and only if for any 
$u, v \in  \la \f S \ra_+$
the vector 
$t \, u + (1-t) \, v \in \la \f S \ra_+$ 
for all 
$0 \le t \le 1 \,.$ 
 On the other hand, by observing that
$u, v \in \la \f S \ra_+ \,,$
we can say that the set  
$\la \f S \ra_+$ 
is convex if and only if for any 
$u, v \in  \la \f S \ra_+$
the vector 
$t \, u + (1-t) \, v \in \la \f S \ra_+$ 
for all 
$0 < t < 1 \,.$

 But, according to the definition of
$\la \f S \ra_+ \,,$ 
we have 
$u = u^1 s_1 + \dots + u^m s_m \,,$
$v = v^1 s'_1 + \dots + v^k s'_k \,,$ 
for $m, k = 1, \dots, n$ 
and
$s_1, \dots, s_m, s'_1, \dots, s'_k \in \f S \,,\,$
$u^1, \dots, u^m, v^1,\dots, v^k \in \Rn^+ \,.$
 Then, 
$tu + (1-t) v = 
\sum_l ((1-t) \, u^l + t \, v^l) \, s_l +
\sum_i (1-t) \, u^i \, s_i + \sum_j t \, v^j \, s'_j \,,$
where
$s_l \in \{s_1, \dots, s_m\} \cap \{s'_1, \dots, s'_k\} \,,$
$s_i \not\in \{s'_1,\dots,s'_k\} \,,$
$s'_j \not\in \{s_1,\dots,s_m\} \,.$
 So, for any 
$0 < t < 1 \,,\,$
$t \, u + (1-t) \, v$ 
is a linear combination of vectors in 
$\{s_1, \dots, s_m\} \cup \{s'_1, \dots, s'_k\} \sub \f S$
with positive coefficients, i.e. it is a vector in 
$\la \f S \ra_+ \,.$\QED
\epf

\bEx\label{cones}
 Let
$\f V$
be a vector space and
$\f S \sub \f V$
a subset consisting of
$n \geq 1$
linearly independent vectors.

 Then, 
$\la \f S \ra_+$
turns out to be a non complete, simple, proper, regular and semi-free
semi--vector space, with the semi--basis 
$\f B = \f S \,.$

 Moreover, any ``transversal section'' of the cone 
$\la\f S \ra_+$
is a polyhedron with
$n$
vertices; they are of the type
$p_i = r^i \, s_i \,,$
with
$r^i \in \Rn^+$
and
$s_i \in \f S \,.$
 Additionally, the boundary of the cone 
$\la\f S \ra_+$
consists of the union of the
$n$
cones generated by the 
$n$
subsets
$\f S_i \byd \f S - \{s_i\} \sub \f S \,,$
with
$s_i \in \f S \,.$

 In particular, 
$\Rn^+$ 
is a non complete, simple, proper, regular and semi-free semi--vector
space with semi--dimension
$1 \,;$
its semi--bases are of the type
$\f B = (s) \,,$
with
$s \in \Rn^+ \,.$
 Clearly, analogous considerations hold for
$\Rn^- \,.$

 More generally,
$(\Rn^+)^n$
is a non complete, simple, proper and regular semi--vector space and 
$(\Rn_0^+)^n$ 
is a complete, simple, proper, regular and
semi--free semi--vector space with semi--dimension
$n \,;$
its semi--bases are of the type
\beq
\f B = 
\big((s^1, 0, \dots, 0) \,,\, \dots \,,\, (0, \dots, 0, s^n)\big)
\sub  (\Rn_0^+)^n \,.
\eeq
 
 Moreover, we have
\beq
\la \f B \ra_+ = (\Rn_0^+)^n - \{0\} \,.
\eeq

 Clearly, analogous considerations hold for
$(\Rn^-)^n \,.$\END
\eEx

\bEx
 Let
$\f V$
be a vector space and
$\f S \sub \f V$
a subset consisting of
$n \geq 1$
linearly independent vectors.

 Then, 
$\la \f S \ra_+$
minus its boundary is a non complete, simple, proper and regular
semi--vector space, but it is not semi--free because it does not
admit any semi--basis.\END
\eEx

\bEx
  Let
$\f V$
be a vector space of dimension
$n$
and
$\f S \sub \f V$
a subset consisting of 
$n$
linearly independent vectors.
 Moreover, let us consider a further non vanishing vector
$v \in \f V - \la \f S \ra_+$
(thus, by hypothesis
$v$
turns out to be a linear combination of
$\f S$).

 Then,
$\la \{v\} \cup \f S \ra_+$
is not semi--free because the uniqueness of the positive decomposition
of semi-vectors is no longer true.

 For instance, a cone in
$\Rn^3 \,,$
whose ``transversal section'' is a square or a disk, turns out to be a
non complete, simple, proper and regular semi--vector space.\END
\eEx

\bEx
 Let
$\f S$
be a set and
$\f U$
a semi--vector space; then, the set
$\Map(\f S,\f U)$
consisting of all maps
$f : \f S \to\f U$
turns out to be a semi--vector space in a natural way.

 Moreover, if
$\f U$
is complete, then
$\Map(\f S,\f U)$
is complete.
 Furthermore, if
$\f U$
is regular, then
$\Map(\f S,\f U)$
is regular.

 Additionally, let 
$\f U$
be complete and semi--free with a semi--basis
$\f B \sub \f U \,.$
 Define the subset
\beq
\f B' \byd \{f_{(\sig,b_i)} \st \sig \in \f S \,,\; b_i \in \f B\}
\sub \Map(\f S,\f U) \,,
\eeq
where, for each
$\sig \in \f S$
and
$b_i \in \f B \,,$
we have defined the map
\beq
f_{(\sig,b_i)} : \f S \to \f U : s \mto \del^\sig_s \, b_i \,,
\ssep{where}
\del^\sig_s = 1
\text{ if } s = \sig
\; \text{ and } \;
\del^\sig_s = 0
\text{ if } s \neq \sig \,,
\eeq

 Then,
$\Map(\f S,\f U)$
is semi--free with a semi--basis 
$\f B' \byd \{f_{(\sig,b_i)}\} \,.$\END
\eEx

\bEx
 Let
$\f V$
be a vector space.
 Then, the set of positive definite symmetric bilinear forms 
$f : \f V \car \f V \to \Rn$
is a proper semi--vector space.\END
\eEx

 Let us exhibit the following non standard example of semi--vector
space.

\bEx\label{union of semi--vector spaces}
 Let us consider two semi--vector spaces
$\f U$
and
$\f V$
and define the set
$\f W \byd \f U \sqcup \f V$
equipped with an addition and a scalar multiplication as follows
\bgt
u_1 + u_2 \,,\,
\lam \, u 
\sep{as in }
\f U \,,
\qquad
u + v = u = v + u \,,
\qquad
v_1 + v_2 \,,\,
\lam \, v 
\sep{as in }
\f V \,,
\end{gather*}
for each
$u, u_1, u_2 \in \f U \,,\, v, v_1, v_2 \in \f V \,,\, 
\lam \in \Rn^+ \,.$

 We can easily see that
$\f W$
turns out to be a semi--vector space.
 In fact, we can verify the axioms step by step.
 For instance, in order to verify the associativity of the addition,
it suffices to consider the following equalities, for each
$u_1, u_2 \in \f U \,,\, v_1, v_2 \in \f V \,,$
\bgt
(u_1 + v_1) + v_2 = u_1 + v_2 = u_1 = u_1 + (v_1 + v_2)
\\
(u_1 + v_1) + u_2 = u_1 + u_2 = u_1 + (v_1 + u_2)
\\
(u_1 + u_2) + v_1 = u_1 + u_2 = u_1 + (u_2 + v_1)
\\
(v_1 + u_1) + u_2 = u_1 + u_2 = v_1 + (u_1 + u_2)
\\
(v_1 + u_1) + v_2 = u_1 + v_2 = u_1 = v_1 + u_1 = v_1 + (u_1 + v_2)
\\
(v_1 + v_2) + u_1 = u_1 = v_1 + u_1 = v_1 + (v_2 + u_1) \,.
\end{gather*}

 We can verify the other axioms in a similar way.

 If
$\f V$
is complete, then
$\f W$
turns out to be complete and we obtain
$0_\f W = 0_\f V \,.$
  
 If both
$\f U$
and
$\f V$
are complete, then
$\f W$
turns out to be complete, but
$0_\f U$
is not a neutral element for 
$\f W \,.$

 If
$\f U$
is complete, then
$\f W$
needs not to be complete.\END
\eEx

\bPr\label{completion of semi--vector spaces}
 If
$\f U$
is a non complete semi--vector space, then we can naturally extend
$\f U$
to a complete semi--vector space
$\haf U$
by considering the complete semi--vector space
$\haf U \byd \f U \sqcup \{0\}$
defined by a procedure as in the above 
Example \ref{union of semi--vector spaces}.\END
\ePr

\bRm
 With reference to the above 
Proposition \ref{completion of semi--vector spaces}, if the
semi--vector space
$\f U$
is already complete, then the original
$0_\f U$
is no longer a neutral element of the new completed semi--vector
space, because we have
$0_\f U + 0 = 0_\f U$
and not
$0_\f U + 0 = 0 \,.$\END
\eRm

\bEx
 With reference to Example \ref{cones}, the sets
$\la \haf S \ra_+ \byd \la\f S \ra_+ \cup \{0\}$
and
$\la \haf S \ra_- \byd \la\f S \ra_- \cup \{0\}$
are complete semi--free and proper semi--vector spaces with
semi--dimension
$n \,.$

 Moreover, 
$\ha\Rn^+ \eqv \Rn^+_0$
and
$\ha\Rn^- \eqv \Rn^-_0$ 
are complete semi--free and proper semi--vector spaces with
semi--dimension
$1 \,.$\END
\eEx

 We can easily define the concepts of semi--vector subspace and
product of semi--vector spaces.
\subsection{Semi--linear maps}
\label{Semi--linear maps}
\bsm
 The notion of semi--linear map is similar to that of linear map.
 However, the lack of some standard properties of semi--vector
spaces yields subtle problems which require additional care.
\esm

\bDf
 A map
$f :\f U \to \f V$
between semi--vector spaces is said to be \emp{semi--linear} if, for
each
$u,v \in\f U, \, r \in \Rn^+ \,,$
we have
$f(u + v) = f(u) + f(v)$
and
$f(r u) = r f(u) \,.$\END
\eDf

 Of course, a linear map between vector spaces is also semi--linear.

 If
$\f U$
and
$\f V$
are semi--vector spaces, then we obtain the semi--vector subspace
\beq
\Slin(\f U,\f V) \byd
\{f :\f U \to \f V \;|\; f \text{ is semi--linear}\} \sub
\Map(\f U,\f V) \,.
\eeq

 In particular, if
$\f U$
is a semi--vector space, then its \emp{semi--dual} is defined to be
the semi--vector space
\beq
\f U^\star \byd \Slin(\f U,\Rn^+) \,.
\eeq

 If
$\f U$
and
$\f V$
are vector spaces,
then
$\Slin (\f U,\f V)$
turns out to be a semi--vector subspace
$\Slin (\f U,\f V) \sub \Lin(\f U,\f V) \,.$

 The composition of two semi--linear maps is semi--linear.
 Hence, semi--vector spaces and semi--linear maps constitute a
category.

 If
$f :\f U \to \f V$
is a bijective
semi--linear map, then the \emp{inverse} map
$f^{-1} :\f V \to \f U$
is also semi--linear.

 If
$f :\f U \to \f V$
is a
semi--linear map, then the \emp{transpose} map
$f^\star : \f V^\star \to \f U^\star : \alp \mto \alp \com f$
is also semi--linear.

\bPr
 Let
$\f U$
be a complete semi--vector space and
$\f V$
a regular and complete semi--vector space.
 If
$f : \f U \to \f V$
is a semi--linear map, then
$f(0_\f U) = 0_\f V \,.$
\ePr

\bpf
 For each
$u \in \f U \,,$
we have
$f(u) + 0_\f V = f(u) = f(u + 0_\f U) = f(u) + f(0_\f U) \,.$
 Hence, we obtain
$f(0_\f U) = 0_\f V \,.$\QED
\epf

\bPr
 Let
$\f U$ be complete semi--vector space
and
$\f V$
be semi--vector space.
 If there exists a surjective semi--linear map 
$f : \f U \to \f V \,,$
then
$\f V$ 
is complete and 
$f(0_\f U) = 0_\f V \,.$
\ePr

\bpf
 In virtue of the surjectivity of
$f \,,$ 
for any 
$v\in \f V$
there exists
$u \in \f U$
such that 
$f(u) = v \,.$
 Then, we have
$v = f(u) = f(u + 0_U) = f(u) + f(0_U) = v + f(0_U)\,,$
which implies that 
$f(0_U) = 0_V \,$ 
is a neutral vector in 
$\f V \,.$\QED
\epf

\bPr
 Let
$\f U$
be a semi--free semi--vector space and
$\f B \sub \f U$
a semi--basis.
 Moreover, let
$\f V$
be any semi--vector space and
$\f S \sub \f V$
any subset.
 Then, there exists a unique semi--linear map
$f : \f U \to \f V$
such that
$f(b_i) = v_i \,,$
with
$b_i \in \f B$
and
$v_i \in \f S \,.$\END
\ePr

\bCr
 Let
$\f U$
be a semi--free semi--vector space of semi--dimension
$n \,.$

 If
$\f U$
is non complete, then it is semi--linearly isomorphic to
$(\Rn^+_0)^n - \{0\} \,.$

 If
$\f U$
is complete, then it is semi--linearly isomorphic to
$(\Rn^+_0)^n \,.$
\eCr

\bpf
 Let 
$\f B = \{b_1, \dots, b_n\}$ 
be a semi--basis of 
$\f U \,.$
 Then the map 
$\f U \to (\Rn_0^+)^n$ 
characterised by 
$b_i \mto (0, \dots,1,\dots ,0) \,,$
with 1 on the 
$i$--th position, defines the semi-linear isomorphism.\QED 
\epf

\bLm\label{semi--linear maps and equality}
 Let
$\f U$
be a semi--free semi--vector space and
$\f V$
a vector space.
 Moreover, let
$u \,,\, \ba u \in \f U \,.$
 If, for any semi--linear map
$f : \f U \to \f V$
we have
$f(u) = f(\ba u) \,,$
then
$u = \ba u \,.$
\eLm

\bpf
 By considering a semi--basis
$\f B \,,$
the equality
$f(u) = f(\ba u)$
implies
$\sum_i u^i \, f(b_i) = \sum_i \ba u^i \, f(b_i) \,.$
 Then, the arbitrariness of
$f$
implies
$u^i = \ba u^i \,.$\QED
\epf

\bPr
 Let
$\f U$
and
$\f U'$
be semi--free semi--vector spaces and let
$\f U'$
be complete; moreover, let
$\f B \sub \f U$
and
$\f B' \sub \f U'$
be semi--bases.
 Then, the semi--vector space
$\Slin(\f U, \f U')$
turns out to be complete and semi--free and with a semi--basis
$\f S \byd \{f_{ij}\} \sub \Slin(\f U, \f U')$
consisting of the semi--linear maps
$f_{ij} : \f U \to \f U'$
uniquely defined by the condition
$f_{ij} : b_h \mto \del_{hi} \, b'_j \,.$\END
\ePr
\subsection{Sesqui--tensor products}
\label{Sesqui--tensor products}
\bsm
 Next, we introduce the notion of tensor product between
semi--vector spaces and vector spaces.
 This construction is quite similar to the analogous construction for
vector spaces, but the lack of standard properties of
semi--vector spaces requires an additional care.
 In particular, it is worth specifying explicitly the different roles
of
$\Rn$
and
$\Rn^+ \,.$
\esm

 Let
$\f U$
be a semi--vector space and
$\f V$
a vector space.

\bDf\label{Definition: sesqui--tensor product}
 A \emp{(right) sesqui--tensor product} between
$\f V$
and
$\f U$
is defined to be a vector space
$\f V \rten \f U$
along with a map
$\rten : \f V \car \f U \to \f V \rten \f U \,,$
which is linear with respect to the 1st factor and semi--linear with
respect to the 2nd factor and which fulfills the following universal
property:

 if
$\f W$
is a vector space and
$f : \f V \car \f U \to \f W$
a map which is linear with respect to the 1st factor and semi--linear
with respect to the 2nd factor, then there exists a unique linear map
$\ti f : \f V \rten \f U \to \f W \,,$
such that
$f = \ti f \com \rten \,.$\END
\eDf

\bTh\label{Theorem: sesqui--tensor product}
 The sesqui--tensor product exists, is unique up to a distinguished
linear isomorphism and is linearly generated by the image of the map
$\rten : \f V \car \f U \to \f V \rten \f U \,.$
\eTh

\bpf
 The proof is analogous to that for the tensor product of vector
spaces, with an additional care.

\emp{Existence.}
 We consider the vector space
$\f F$
consisting of all maps
$\phi : \f V \car \f U \to \Rn \,,$
which vanish everywhere except on a finite subset of
$\f V \car \f U \,.$
 Clearly, the set
$\f F$
becomes a vector space in a natural way.
 Accordingly, each
$\phi \in \f F$
can be written as a formal sum of the type
\beq
\phi = \phi^{11} (v_1,u_1) + \dots + \phi^{nm}  (v_n,u_m) \,,
\eeq
where
$\phi^{ij} \eqv \phi(v_i,u_j) \in \Rn$
are the (possibly) non vanishing values of
$\phi \,.$

 Next, we consider the subset
$\f S \sub \f F$
consisting of elements of the type
\bgt
(v + v', \, u) - (v,u) - (v',u) \,,
\qquad
(v, \, u + u') - (v,u) - (v,u') \,,
\\
(sv, u) - s(v, u) \,,
\qquad
(v, ru) - r(v, u) \,,
\end{gather*}
with
$v,v' \in \f V \,,$
$u, u' \in \f U \,,$
$s \in \Rn \,,$
$r \in \Rn^+ \,.$
 Then, we consider the vector subspace
$\la \f S \ra_{\Rn} \sub \f F$
linearly generated by
$\f S$
on
$\Rn \,.$
 Eventually, we obtain the quotient vector space and the bilinear map
\beq
\f V \rten \f U \eqv
\f F/ \la \f S \ra_{\Rn}
\ssep{and}
\rten \byd q \com \jmath  : \f V \car \f U \to \f V \rten \f U \,,
\eeq
where
$\jmath : \f V \car \f U  \hto \f F$
and
$q : \f F \to 
\f F/ \la \f S\ra_{\Rn}$
are the natural inclusion and the quotient projection.

 Clearly,
$\f V \rten \f U$
is linearly generated by the image of the map
$\rten \,.$

 Now, let us refer to the universal property (Definition
\ref{Definition: sesqui--tensor product}). 
 If the linear map
$\ti f : \f V \rten \f U \to \f W$
such that
$f = \ti f \com \rten$
exists, then it is unique because
$\f V \rten \f U$
is linearly generated by the image of the map
$\rten \,.$
 Indeed, such a map exists. 
 In fact, we can easily prove that the bilinear map
$f : \f V \car \f U \to \f W$
yields naturally a linear map
$f' : \f F \to \f W \,,$
which passes to the quotient yielding the required linear map
$\ti f : \f V \rten \f U \to \f W \,.$

\myskip

\emp{Uniqueness.}
 The sesqui--tensor product is ``unique'' in the following sense.
If
$\f V \rten \f U$
and
$\f V \bau\rten \f U$
are sesqui--tensor products, then the universal properties of the two
sesqui--tensor products yield the following commutative diagram
\bdg[size=2em]
\f V \rten \f U
&&\rLR^{\wti{\bau\rten}}_{\wti{\rten}}
&&\f V \bau\rten \f U
\\
&\luTo_\rten
&&\ruTo_{\bau\rten}
\\
&&\f V \car \f U
\edg
where
$\wti{\bau\rten} :
\f V \rten \f U \to \f V \bau\rten \f U$
and
$\wti\rten :
\f V \bau\rten \f U \to \f V \rten \f U$
are mutually inverse linear isomorphisms.\QED
\epf

\bNt
 Clearly, for each
$v,v' \in \f V \,,\,$
$u, u' \in \f U \,,\,$
$s \in \Rn \,,\,$
$r \in \Rn^+ \,,$
we have
\bgt
(v + v') \rten u = v \rten u + v' \rten u \,,
\qquad
v \rten (u + u') = v \rten u + v \rten u' \,,
\quad
\\
(s \, v) \rten u = s \, (v \rten u) \,,
\qquad
v \rten (r \, u) = r \, (v \rten u) \,.\END
\end{gather*}
\eNt

 The following annihilation rules follow from the universal property:

\bPr
 If
$u \in \f U \,,\, 0_\f V \in \f V$
then we have
$0_\f V \rten u = 0 \in \f V \rten \f U \,.$
\ePr

\bpf
 If
$f : \f V \rten \f U \to \f W$
is any linear map, then
$\phi \byd f \com \rten : \f V \car \f U \to \f W$
is a linear map with respect to the 1st factor.
 Hence, we have
$\phi (0_\f V, u) = 0 \,,$
which, in virtue of the universal property, implies
$f (0_\f V \rten u) = 0_\f W \,.$
 Therefore, in virtue of the arbitrariness of
$f \,,$
we obtain
$0_\f V \rten u = 0 \,.$\QED
\epf

\bPr
 Let
$\f U$
be complete.
 If
$0_\f U \in \f U \,,\, v \in \f V \,,$
then we have
$v \rten 0_\f U = 0 \in \f V \rten \f U \,.$
\ePr

\bpf
 If
$f : \f V \rten \f U \to \f W$
is any linear map, then
$\phi \byd f \com \rten : \f V \car \f U \to \f W$
is a semi--linear map with respect to the 2nd factor.
 Hence, we have
$\phi (v, 0_\f U) = 0_\f W \,,$
which, in virtue of the universal property, implies
$f (v \rten 0_\f U) = 0_\f W \,.$
 Therefore, in virtue of the arbitrariness of the linear map
$f \,,$
we obtain
$v \rten 0_\f U = 0 \,.$\QED
\epf

\bPr
 Let
$\f U$
be semi-free.
 Moreover, let
$\f B$
be a basis of
$\f V$
and
$\f C$
a semi--basis of
$\f U \,.$
 Then,
\beq
\f B \rten \f C \byd \{b_i \rten c_j \st 
b_i \in \f B \,,\, c_j \in \f C\}
\eeq
is a basis of
$\f V \rten \f U \,.$
 Thus, we have
\beq
\dim (\f V \rten \f U) = (\dim \f V) \, (\sdim \f U) \,.
\eeq
\ePr

\bpf
 Clearly,
$\f B \rten \f C$
linearly generates
$\f V \rten \f U \,.$

 Next, let us prove that the elements of
$\f B \rten \f C$
are linearly independent.
 For this purpose, let us observe that the universal property of the
sesqui--linear tensor product yields a bijection 
$f \mto \ti f$
between the maps
$f : \f V \car \f U \to \Rn$
which are linear with respect to the 1st factor and semi--linear with
respect to the 2nd factor and the linear maps
$\ti f : \f V \rten \f U \to \Rn \,,$
according to the rule
$\ti f(v \rten u) = f (v, u) \,,$
for each
$v \in \f V \,,\, u \in \f U \,.$
 Now, let us consider an element
$t \byd \sum_{ij} t^{ij} b_i \rten c_j \in \f V \rten \f U \,.$
 Indeed, for any
$\ti f$
as above, we have
$\ti f(\sum_{ij} t^{ij} b_i \rten c_j) = 
\sum_{ij} t^{ij} f(b_i, c_j) \,.$
 Then, a property of semi--bases and bases implies that
$\sum_{ij} t^{ij} f(b_i, c_j) = 0 \,,$
for all
$f$
as above, if and only if
$t^{ij} = 0 \,.$
 Hence,
$\ti f(\sum_{ij} t^{ij} b_i \rten c_j) = 0 \,,$
for all
$\ti f$
as above, if and only if
$t^{ij} = 0 \,.$
 Therefore, a property of bases implies that
$(b_i \rten c_j)$
is a basis of the sesqui--tensor product.\QED
\epf

\bRm
 If
$\f U$
and
$\f V$
are vector spaces, then we can consider their sesqui--tensor product
and tensor product,  by considering, respectively, one factor as
a semi--vector spaces, or both factors as vector spaces.
 We stress that the above tensor products are different.
 For instance, if
$v \in \f V$
and
$u \in \f U \,,$
then
\bat{6}
&- 
(v \rten u) 
&&\neq 
v \rten (-u) 
&&\neq 
(-v) \rten u 
&&= 
- (v \rten u) \,,
&&\qquad
(-v) \rten (-u) 
&&\neq v \rten u \,,
\\
&- 
(v \ten u) 
&&= 
v \ten (-u) 
&&= (-v) \ten u 
&&= - (v \ten u) \,,
&&\qquad
(-v) \ten (-u) 
&&= 
v \ten u \,.\END
\end{alignat*}
\eRm

\bPr
 Let
$\f V$
and
$\f U$
be vector spaces. 
 Then, we have
\beq
\dim (\f V \rten \f U) = 2 \, (\dim \f V) \, (\dim \f U) \,.
\eeq

 Indeed, the universal properties of the tensor products yield the
natural surjective semi--linear map
\beq
\gr\pi : \f V \rten \f U \to \f V \ten \f U :
v \rten u \mto v \ten u \,,
\eeq
whose kernel is semi--linearly generated by the elements of the type
\beq
v \rten u + v \rten (-u) \in \f V \rten \f U \,,
\eeq
with
$v \in \f V \,,\, u \in \f U \,.$
 Thus, we have
$\dim \, (\ker \gr\pi) = (\dim \f V) \, (\dim \f U) \,.$\END
\ePr

\myskip

 We can introduce the left sesqui--tensor product
$\f U \lten \f V$
analogously to the right sesqui--tensor product
$\f V \rten \f U \,.$
 Clearly, we have a natural linear isomorphism
$\f V \rten \f U \seq \f U \lten \f V \,,$
which is characterised by the map
$v \rten u \mto u \lten v \,.$
\subsection{Universal vectorialising space}
\label{Universal vectorialising space}
\bsm
 The sesqui--tensor product of a semi--free semi--vector space with
$\Rn$
yields a distinguished extension of the semi--vector space into
a vector space.
\esm

 Let us consider a semi--vector space
$\f U$
and the distinguished vector space
$\Rn \rten \f U \,.$

\bLm
 Each element 
$t \in \Rn \rten \f U$
can be written as
\[
t = 
a \, 1 \rten u_+ + b \, (-1) \rten u_- = 
a \, 1 \rten u_+ - b \, 1 \rten u_-\,,
\]
where each of the the parameters
$a$
and
$b$
might assume the value
$1$
or
$0$
and
$u_+ \,, u_- \in \f U \,.$
\eLm

\bpf
 Each element
$t \in \Rn \rten \f U$
can be written as
\beq
t = \sum_{ij} t^{ij} (\lam_i \rten u_j) = 
\sum_j (\sum_i t^{ij} \lam_i) \rten u_j =
\sum_j r^j \rten u_j \,,
\eeq
with
$t^{ij}, \lam_i \in \Rn \,,$
$u_j \in \f U$
and
$r^j \byd \sum_i t^{ij} \lam_i \in \Rn \,.$

 Hence, we obtain the result by recalling that
$0 \rten u = 0$
and by setting
\beq
u_+ \byd \sum_j r^j \, u_j \,,
\sep{for}
r^j > 0
\ssep{and}
u_- \byd \sum_j (- r^j) \, u_j \,,
\sep{for}
r^j < 0 \,.\QED
\eeq
\epf

\bPr
 The sesqui--tensor product
$\Rn \rten \f U$
is the union of three semi--vector subspaces
\beq
\Rn \rten \f U = \f U_+ \cup \{0\} \cup \f U_- \,,
\eeq
where
\beq
\f U_+ \byd \{1 \rten u \st u \in \f U\} \,,
\qquad
\f U_- \byd \{- 1 \rten u \st u \in \f U\} \,.
\eeq

 Thus,
$\Rn \rten \f U$
is linearly generated by the subset consisting of elements of the
type
$1 \rten u \,,$
with 
$u \in \f U \,.$

 Moreover, we have the natural semi--linear map
\beq
\imath : \f U \to \Rn \rten \f U : u \mto 1 \rten u \,.\END
\eeq
\ePr

\bLm
 Let
$\f U$
be semi--free.
 Then, for each non neutral element
$u \in \f U \,,$
we have
$1 \rten u \neq 0$
and
$(-1) \rten u \neq 0 \,.$
\eLm

\bpf
 By the hypothesis that
$\f U$
is semi--free, there exist a semi--linear map 
$\phi : \f U \to \Rn$
such that
$\phi(u) \neq 0 \,.$
 This map yields the map
$f : \Rn \car \f U \to \Rn : (r, v) \mto r \, \phi(v) \,;$
indeed, the map 
$f$
is linear with respect to the 1st factor and semi--linear with
respect to the 2nd factor and fulfills the equalities
$f(1,u) = \phi(u)$
and
$f(-1,u) = - \phi(u) \,.$
 Then, in virtue of the universal property, we obtain the linear map
$\wti f : \Rn \rten \f U \to \Rn \,.$
 Thus, we have
$\wti f (1 \rten u) = \phi (u) \neq 0$
and
$\wti f ((-1) \rten u) = -\phi (u) \neq 0 \,;$
hence
$1 \rten u \neq 0$
and
$(-1) \rten u \neq 0 \,.$\QED
\epf

\bLm
 Let
$\f U$
be semi--free.
 Then, for each non neutral elements
$u \,, v \in \f U \,,$
we have
$1 \rten u \neq (-1) \rten v \,.$
\eLm

\bpf
 By the hypothesis that
$\f U$
is semi--free, there exists a semi--linear map
$\phi : \f U \to \Rn \,,$
with positive values, excepts on the possible neutral element of
$\f U \,.$
 This map yields the map
$f : \Rn \car \f U \to \Rn : (r, v) \mto r \, \phi(v) \,;$
indeed, the map 
$f$
is linear with respect to the 1st factor and semi--linear with
respect to the 2nd factor.
 Then, in virtue of the universal property, we obtain the linear map
$\wti f : \Rn \rten \f U \to \Rn \,.$
 Indeed, we obtain
$\ti f(1 \rten u) = \phi (u) \in \Rn^+$
and
$\ti f((-1) \rten v) = - \phi (v) \in \Rn^- \,.$
 Hence,
$1 \rten u \neq (-1) \rten v \,,$
because there is a linear map which takes different values on these
elements.\QED
\epf

\bLm
 Let
$\f U$
be semi--free.
 Then, for each
$u \,, \ba u \in \f U \,,$
$1 \rten u = 1 \rten \ba u$
implies
$u = \ba u \,.$
\eLm

\bpf
 Let us consider any semi--linear map
$\phi : \f U \to \Rn$
and the induced map
$f : \Rn \car \f U \to \Rn : (r, u) \mto r \, \phi(u) \,,$
which is linear with respect to the 1st factor and semi--linear
with respect to the 2nd factor.
 Thus, we obtain  the induced linear map
$\ti f : \Rn \rten \f U \to \Rn : r \rten u \mto r \, \phi(u) \,.$

 Then, the equality
$1 \rten u = 1 \rten \ba u$
implies
$\ti f(1 \rten u) = \ti f (1 \rten \ba u) \,,$
hence
$\phi(u) = \phi(\ba u) \,.$
 Therefore, the arbitrariness of
$\phi$
implies
$u = \ba u \,,$
in virtue of Lemma \ref{semi--linear maps and equality}.\QED
\epf

\bNt
 Let
$\f U$
be semi--free.

 If
$\f U$
is non complete, then
$\f U_+ \cap \f U_- = \emptyset \,.$

 If
$\f U$
is complete, then
$\f U_+ \cap \f U_- = \{0\} \,.$\END
\eNt

\bPr\label{Proposition: universal vectorialising space}
 Let
$\f U$
be semi--free.
 Then, the map
$\imath$
is injective.
 Moreover, if
$\f B$
is a semi--basis of
$\f U \,,$
then
$(1 \rten b_i)$
is a basis of
$\Rn \rten \f U \,.$
 Hence,
\beq
\dim (\Rn \rten \f U) = \sdim \f U \,.\END
\eeq
\ePr

\bDf
 If
$\f U$
is semi--free, then we say
$\baf U \byd \Rn \rten \f U$
to be the
\emp{universal vector extension} of
$\f U \,.$\END
\eDf

 In fact, this vector space allows us to transform any semi--linear
map into a linear map.
 Moreover, this space is the smallest one with this property.

\bPr
Let
$\f U$
be semi--free.
 Moreover, let
$\f W$
be a vector space and
$f : \f U \to \f W$
a semi--linear map.
 Then, the following facts hold.

1) There exists a unique linear map
$\ba f : \baf U \to \f W \,,$
such that
$f = \ba f \com \imath \,.$
 This map is given by
\beq
\ba f(1 \rten u) = f(u) \,,
\quad
\ba f ((-1) \rten u) = - f(u) \,,
\quad
\Al u \in \f U \,.
\eeq

2) If
$\baf U'$
is another vector space and
$\imath' : \f U \hto \baf U'$
a semi--linear inclusion which fulfill the above property, then there
exist a unique linear inclusion
$\rho : \baf U \hto \baf U' \,,$
such that
$\imath' = \rho \com \imath$
and
$\ba f = \ba f' \com \rho \,.$
\ePr

\bpf
 1) Let us consider the map
$\phi : \Rn \car \f U \to \f W : (r, u) \mto r \, f(u) \,,$
which is linear with respect to the 1st factor and semi--linear with
respect to the 2nd factor.
 Then, in virtue of the universal property of the sesqui--tensor
product, there is a unique linear map
$\ba f : \baf U \to \f W \,,$
such that
$\ba f (r \rten u) = \phi (r, u) \byd r \, f (u) \,.$
 In particular, we have
$\ba f(1 \rten u) = f(u) \,.$

 Moreover, if
$\ba f' : \baf U \to \f W$
is another linear map such that
$\ba f'(1 \rten u) = f(u) \,,$
then
$\ba f'(1 \rten u) = f(u) = \ba f (1 \rten u)$
implies also
$\ba f'((-1) \rten u) = - f(u) = \ba f ((-1) \rten u) \,.$
 Hence, being
$\Rn \rten \f U$
semi--linearly generated by the elements of the type
$1 \rten u$
and
$(-1) \rten u \,,$
we obtain
$\ba f' = \ba f \,.$

 2) The existence and uniqueness of the semi--linear map
$\rho : \baf U \to \baf U'$
follow from the universal property of the sesqui--tensor product
$\Rn \rten \f U \,,$
by considering the map
$\Rn \car \f U \to \baf U' : (r, u) \mto r \, \imath'(u) \,.$\QED
\epf

 For instance, the sesqui--tensor product
$\baf U' \byd \Cn \rten \f U$
is also a vector extension of
$\f U$
and we have the distinguished real linear inclusion
$\Rn \rten \f U \hto \Cn \rten \f U \,.$

 In an analogous way we can prove the following fact.

\bPr
 Let
$\f U$
and
$\f V$
be semi-free semi--vector spaces.
 Moreover, let
$\f W$
be a vector space and
$f : \f U \car \f V \to \f W$
a semi--bilinear map.
 Then, the following facts hold:

 1) There exists a unique bilinear map
$\ba f : \baf U \car \baf V \to \f W \,,$
such that
$f = \ba f \com (\imath_\f U \car \imath_\f V) \,.$

 This map is given by
\bgt
\ba f(1 \rten u \,,\, 1 \rten v) = f(u, v) \,,
\quad
\ba f ((-1) \rten u \,,\, 1 \rten v) =
\ba f (1 \rten u \,,\, (-1) \rten v) = - f(u, v) \,,
\\
\ba f ((-1) \rten u \,,\, (-1) \rten v) = f(u, v) \,,
\qquad
\Al (u, v) \in \f U \car \f V \,.
\end{gather*}

2) If
$\baf U' \car \baf V'$
is another vector space and
$(\imath'_\f U \car \imath'_\f V) : 
\f U \car \f V  \hto \baf U' \car \baf V'$
a semi--linear inclusion which fulfill the above property, then there
exist a unique linear inclusion
$\rho : \baf U \car \baf V \hto \baf U' \car \baf V'\,,$
such that
$(\imath'_\f U \car \imath'_\f V) = 
\rho \com (\imath_\f U \car \imath_\f V)$
and
$\ba f = \ba f' \com \rho \,.$\END
\ePr

\bNt
 From Proposition \ref{Proposition: universal vectorialising space} 
and 
Example \ref{cones} it follows that all semi--free semi--vector spaces
can be regarded as cones in a vector space.\END
\eNt

\bRm
 Let
$\f V$
be a vector space.
 Then the map
$\imath : \f V \to \Rn \rten \f V$
is not a linear isomorphism; in fact, for each
$v \in \f V \,,$
we have
\beq
\imath (- v) = 1 \rten (- v) \neq - (1 \rten v) = -\imath (v) \,.
\eeq

 Indeed, if
$\f B = (b_1, \dots b_n)$
is a basis of
$\f V \,,$
then
\beq
\big(1 \rten b_1, \dots, 1\rten b_n, \;
(-1) \rten b_1, \dots, (-1) \rten b_n\big) \sub \Rn \rten \f V
\eeq
is a basis of
$\Rn \rten \f V \,.$
 
 Hence, we have
\beq
\dim (\Rn \rten \f V) = 2 \, \dim \f V \,.
\eeq

 Thus, we have
\beq
\baf V \byd \Rn \rten \f V \neq \f V \seq \Rn \ten \f V \,.
\eeq

 On the other hand, we have a surjective linear map
\beq
\baf V \to \f V :
r \rten v \mto r \, v \,,
\eeq
whose kernel is linearly generated by the elements of the type
\beq
1 \rten v + 1 \rten (- v) \,,
\ssep{with}
v \in \f V \,.\END
\eeq
\eRm

\bPr
 Let
$\f U$
and
$\f V$
be semi--free semi--vector spaces.
 Then, we have the natural injective linear map
\beq
\ol{\Slin(\f U, \f V)} \to \Lin(\baf U, \baf V) : 
r \rten f \mto \phi \,,
\sep{where}
\phi : \baf U \to \baf V : s \rten u \mto r \, (s \rten f(u)) \,.
\eeq

 Indeed, we have
\bal
\dim \, (\ol{\Slin(\f U, \f V)}) 
&= 
\sdim \, (\f U) \,  \, \sdim \, (\f V) \,,
\\
\dim \, (\Lin(\baf U, \baf V)) 
&= 
\sdim \, (\f U) \,  \, \sdim \, (\f V) \,.\END
\end{align*}
\ePr

\bCr
 Let
$\f U$
be a semi--free semi--vector space.
 Then, we have a natural injective linear map
\beq
\ol{\f U^\star} \to \baf U^* \,.\END
\eeq
\eCr
\subsection{Semi--tensor products}
\label{Semi--tensor product}
\bsm
 Now, we are in a position to introduce the notion of tensor
product between semi--free semi--vector spaces.

 The procedure is similar to that for vector spaces, but we need to
pass through the universal vector extension, in order to overcome
the lack of standard properties of semi--vector spaces.
\esm

 Let
$\f U$
and
$\f V$
be semi--free semi--vector spaces.

\bDf\label{Definition: semi--tensor product}
 A \emp{semi--tensor product} between
$\f U$
and
$\f V$
is defined to be a semi--vector space
$\f U \sten \f V$
along with a semi--bilinear map
$\sten : \f U \car \f V \to \f U \sten \f V \,,$
which fulfills the following universal property:

 if
$\f W$
is a semi--vector space and
$f : \f U \car \f V \to \f W$
a semi--bilinear
map, then there exists a unique semi--linear map
$\ti f :\f U \sten \f V \to \f W \,,$
such that
$f = \ti f \com \sten  \,.$\END
\eDf

\bTh\label{Theorem: semitensor product}
 The semi--tensor product exists and is unique up to a distinguished
semi--linear isomorphism.
\eTh

\bpf
 The uniqueness can be proved by a standard procedure as in Theorem
\ref{Theorem: sesqui--tensor product}.
 Then, we have to prove the existence of the semi--tensor product.

 For this purpose, we consider the subset
$\f U \sten \f V \sub \baf U \rten \f V$
consisting of the semi--linear combinations of elements of the type
$(1 \rten u) \rten v \,,$
with
$u \in \f U$
and
$v \in \f V \,,$
and the map
$\sten : \f U \car \f V \to \f U \sten \f V : 
(u, v) \mto (1 \rten u) \rten v \,.$

 We can easily see that
$\f U \sten \f V$
is a semi--vector space and
$\sten : \f U \car \f V \to\f U \sten \f V$
a semi--bilinear map.

 Next, we prove that the above objects fulfill the required universal
property.

 Clearly, if the map
$\ti f : \f U \sten \f V \to \f W$
of the universal property exists, then it is unique because
$\f U \sten \f V$
is semi--linearly generated by the image of the map
$\sten \,.$

 Moreover, this map is well defined by the equality
$\ti f((1 \rten u) \rten v) = f(u,v) \,,$
according to the following commutative
diagram
\bdg[size=2em]
\f U \car \f V
&&&\rTo^i
&&&{\baf U \car \baf V}
\\
&\rdTo^f
&&&&\ldTo^{\ol{\imath \com f}}
\\
\dTo^\sten
&&\f W
&\rTo^\imath
&{\baf W}
&&\dTo_{\ten}
\\
&\ruTo_{\ti f}
&&&&\luTo_{\wti{\ol{\imath \com f}}}
\\
\f U \sten \f V
&&&\rTo_j
&&&{\baf U \ten \baf V}
\edg
where the maps
$i, \imath, j$
are natural inclusions and where the map
$(\wti{\ol{\imath \com f}}) \com j$
uniquely factorises through a map
$\ti f \,.$
We can easily see that the map
$\ti f$ 
is
semi--linear and that
$\ti f \com \sten = f \,.$\QED
\epf

 In an analogous way, we can construct the semi--tensor product via
the left sesqui-tensor product (instead of via the right
sesqui--tensor product).
 We can also easily prove that the two constructions are naturally
isomorphic.

\bPr
 The semi--tensor product is a semi--free semi--vector space.
 Moreover, if
$\f B$
and
$\f C$
are semi--bases of
$\f U$
and
$\f V \,,$
respectively, then
\beq
\f B \sten \f C \byd \{b_i \sten c_j \st 
b_i \in \f B \,,\, c_j \in \f C\}
\eeq
is a semi--basis of
$\f U \sten \f V \,.$
 Moreover, we have
\beq
\sdim (\f U \sten \f V) = \sdim \f U \cdot \sdim \f V \,.\END
\eeq
\ePr

\bNt
  We have the natural semi--linear isomorphisms
\beq
\Rn^+ \sten \f U \to \f U : r \sten u \mto r u
\ssep{and}
\f U \sten \Rn^+ \to \f U : u \sten r \mto r u \,.\END
\eeq
\eNt

\bNt
 We have a distinguished semi--linear isomorphism
\beq
\ol{\f U \sten \f V}
\;\seq\;
\baf U \ten \baf V \,.\END
\eeq
\eNt

\bPr
 We have the natural semi--linear inclusion
\beq
\f U^* \sten \f V \hto \Slin(\f U,\f V) \,,
\eeq
characterised by
\beq
\alp \sten v : \f U \to \f V : u \mto \alp(u) \, v \,,
\qquad
\Al \alp \in \f U^*, \, v \in \f V \,.
\eeq

 Moreover, if the semi--dimensions of
$\f U$
and
$\f V$
are finite, then the above inclusion is a semi--linear isomorphism.
\ePr

\bpf
 It follows easily from the universal property in the standard
way.\QED
\epf

 We obtain the contravariant ``semi--tensor algebra'' of a
semi--free semi--vector space in a way analogous to that of vector
spaces.

 Let
$m$
be a positive integer.
 If
$\f U_1, \dots, \f U_m$
are semi--free semi--vector spaces, then we can easily define
the semi--tensor product
$\f U_1 \sten \dots \sten \f U_m$
and prove its basic properties along the same lines as for vector
spaces.
 In particular, if
$\f U = \f U_1 = \dots = \f U_m \,,$
then we set
$\sten^m\f U \byd \f U_1 \sten \dots \sten \f U_m$
and
$\sten^0\f U \byd \Rn^+ \,.$
 Moreover, the semi--direct sum
$\bigoplus_{m \in \B N} \f U^m$
turns out to be a ``semi--algebra'' over
$\Rn^+ \,.$
\section{Positive spaces}
\label{Positive spaces}
\bsm
 The positive spaces constitute a distinguished elementary type of
semi--vector spaces, which allows us to introduce further notions,
such as rational maps and rational powers.

 We use these spaces for achieving an algebraic model of scales and
units of measurement.
 In fact, this is one of the main goals of this paper.
\esm
\subsection{The notion of positive space}
\label{The notion of positive space}
\bsm
 We start by introducing the positive spaces and their basic
properties.
\esm

\bDf
 A \emp{positive space} is defined to be a non--complete semi--free
semi-vector space
$\B U$
of dimension 1.\END
\eDf

\bNt
 In other words, a positive space is a semi-vector space generated
over
$\Rn^+$
by a non vanishing element.

 Thus, each positive space
$\B U$
is semi--linearly isomorphic to
$\Rn^+ \,.$
More precisely, each semi--linear isomorphism
$\B U \to \Rn^+$
is of the type
$\B U \to \Rn^+ : u \mto u/b \,,$
where
$b \in \B U$
and
$r \byd u/b \in \Rn^+$
is the unique positive number such that
$r b = u \,.$

 Clearly, for each
$b \in \B U \,,$
the subset
$(b) \sub \B U$
turns out to be a semi--basis.

 If
$\B U$
and
$\B V$
are positive spaces, then each semi--linear map
$f :\B U \to \B V$
is an isomorphism.
 In fact, we can easily see that
$f$
is surjective and injective.\END
\eNt

 Positive spaces and semi--linear maps constitute a category.

\bNt
 Let
$\B U$
be a positive space.
 The scalar multiplication
$s : \Rn^+ \car \B U \to \B U$
turns out to be a free and transitive action of the group
$(\Rn^+ \,, \cdot)$
on the set
$\B U \,.$\END
\eNt

\bNt
 The semi--tensor product of positive spaces is a positive space.

 In particular, if
$\B U$
is a positive space, then
$\sten^n \B U$
is a positive space.
 Moreover, each element
$t \in \sten^n \B U$
is decomposable; even more, it can be uniquely written as
$t = u \sten \dots \sten u \,,$
with
$u \in \B U \,.$\END
\eNt

 For positive spaces we shall often adopt a notation similar to the
standard notation used for numbers.
 Namely, if
$\B U$
and
$\B U'$ 
are positive spaces, we shall often write
$u \, u' \eqv u \sten u' \in \B U \sten \B U'\,,$
for each
$u \in \B U$
and
$u' \in \B U' \,.$

 Moreover, if
$\B U$
is a positive space and
$u \in \B U \,,$
then the unique element
$1/u \in \B U^\star \,,$
such that
$\la 1/u, \, u \ra = 1 \,,$
is called the \emp{inverse} of
$u$
(not to be confused with the additive inverse).
 Clearly, for each
$u \in \B U$
and
$r \in \Rn^+ \,,$
we have
$
\tfr1{ru} = \tfr1r \tfr1u \,.$
 Moreover,
$(1/u)$
is just the dual semi--basis of
$(u) \,.$
\subsection{Rational maps between positive spaces}
\label{Rational maps between positive spaces}
\bsm
 Next, we discuss the notion of $q$--rational maps between positive
spaces.
\esm

 Let us consider two positive space
$\B U$
and
$\B V$
and a rational number
$q \in \B Q \,.$

\myskip

\bDf
 A map
$f : \B U \to \B V$
is said to be \emp{$q$--rational} (or, \emp{rational of degree}
$q$)
if, for each
$u \in \B U$
and
$r \in \Rn^+ \,,$
we have
$f(r \, u) = r^q \, f(u) \,.$\END
\eDf

 We denote by
$\Rat^q (\B U, \B V) \sub \Map(\B U, \B V)$
the subspace of
$q$--rational maps
between the positive spaces
$\B U$
and
$\B V \,.$

\bPr
 If
$u \in \B U$
and
$v \in \B V \,,$
then there exists a unique
$q$--rational map
$f : \B U \to \B V \,,$
such that
$f(u) = v \,.$\END
\ePr

 The composition of two rational maps is a rational map, whose degree
is the product of the degrees.
 Hence, positive spaces and rational maps constitute a category.

\bNt
 Let 
$q'$
be another rational number.
 If
$f : \B U \to \Rn^+$
is a 
$q$--rational map, then the map
$f^{q'} : \B U \to \Rn^+ : u \mto (f(u))^{q'}$
is 
$(q q')$--rational.\END
\eNt

\bPr
 The subspace
$\Rat^q (\B U, \B V) \sub \Map(\B U, \B V)$
turns out to be a semi-vector subspace and a positive space.
\ePr

\bpf
 The 1st statement is trivial.
 Moreover,
$\Rat^q (\B U, \B V)$
is a positive space because,
for any given
$u \in \B U \,,$
the map
$\Rat^q (\B U, \B V) \to \B V : f \mto f(u)$
is a semi--linear isomorphism.\QED
\epf

\bCr
 A $q$--rational map
$f : \B U \to \B V$
is a bijection if and only if
$q \neq 0 \,;$
in this case the inverse map is
$(1/q)$--rational.\END
\eCr

\bEx
 We have the following distinguished cases.

a) The 0--rational maps
$f : \B U \to \B V$
are just the constant maps.
 Hence, we have the natural semi--linear isomorphism
\beq
\Rat^0 (\B U, \B V) \seq \B V : f \mto f(u) \,,
\eeq
which turns out to be independent of the choice of
$u \in \B U \,.$ 
 In particular, we have
$\Rat^0(\B U, \Rn^+) \seq \Rn^+ \,.$

b) The 1--rational maps
$f : \B U \to \B V$
are just the semi--linear maps.
 Hence, we can write
\beq
\Rat^1 (\B U, \B V) = \Slin (\B U, \B V) \,.
\eeq

 In particular, we have
$\Rat^1(\B U, \Rn^+) = \Slin(\B U, \Rn^+) = \B U^\star.$

c) The $(-1)$--rational maps
$f : \B U \to \B V$
can be identified with the semi--linear maps
$\bau f : \B U^\star \to \B V \,,$
through the natural semi--linear isomorphism
\beq
\Rat^{-1} (\B U, \B V) \to \Slin (\B U^\star, \B V) : 
f \mto \bau f \,,
\eeq
where
$\bau f : \B U^\star \to \B V$
is the unique semi--linear map such that
$\bau f(1/u) = f(u) \,,$
with reference to a chosen element
$u \in \B U \,.$
 Indeed, this isomorphism turns out to be independent of the choice of
$u \in \B U \,.$\END

 In particular, the map
\beq
\inv : \B U \to \B U^\star : u \mto 1/u \,,
\eeq
which associates with each element
$u \in \B U$
its dual form
$1/u \in \B U^\star \,,$
is a $(-1)$--rational map.

 Indeed,
$\inv \in \Rat^{-1} (\B U, \B U^\star)$
is the distinguished element which corresponds to the element
$\id_{\B U^\star} \in \Slin(\B U^\star, \B U^\star) \,,$
through the isomorphism
$\Rat^{-1} (\B U, \B U) \seq \Slin (\B U^\star, \B U) \,.$

 We have also the map
\beq
\inv : \B U^\star \to \B U^\star{}^\star \seq \B U \,.\END
\eeq
\eEx
\subsection{Rational powers of a positive space}
\label{Rational powers of a positive space}
\bsm
 Eventually, we introduce the rational powers of a positive space.

 The basic idea is quite simple and could be achieved in an elementary
way, by referring to a base and showing that the result is independent
of this choice.

 However, a full understanding of this concept is more subtle than it
might appear at first insight and suggests a more sophisticated formal
approach.
\esm

 Let us consider a positive space
$\B U$
and a rational number
$q \in \B Q \,.$

\myskip

\bLm
 The map
\beq
\pi^q : \B U \to \Rat^q (\B U^\star, \Rn^+) : u \mto u^q \,,
\eeq
where
$u^q \in \Rat^q (\B U^\star, \Rn^+)$
is the unique element such that
$u^q (1/u) = 1 \,,$
turns out to be
$q$--rational.
\eLm

\bpf
 In fact, we have
$1 = u^q (1/u)$
and
$1 =
(ru)^q (1/(ru)) =
(ru)^q (1/r \, 1/u) =
(1/r)^q \, (ru)^q (1/u) \,.$

 Hence, we obtain
$u^q (1/u) = (1/r)^q \, (ru)^q (1/u) \,,$
which yields
$(ru)^q = r^q u^q \,.$\QED
\epf

\bDf
 The \emp{$q$--power} of
$\B U$
is defined to be the pair
$(\B U^q, \, \pi^q)$
defined by
\beq
\B U^q \byd \Rat^q (\B U^\star, \Rn^+)
\ssep{and}
\pi^q : \B U \to \B U^q : u \mto u^q \,,
\eeq
where
$u^q : \B U^\star \to \Rn^+$
is the unique
$q$--rational map such that
$u^q (1/u) = 1 \,.$\END
\eDf

\myskip

 We can re--interpret the above notion in a natural way in terms of
semi--tensor powers as follows.

\bNt
 Clearly, for each
$r \in \Rn^+ \,,$
the following diagram commutes
\bcd
\B U @>>> \B U^q
\\
@V{s_r}VV   @VV{s_{r^q}}V
\\
\B U   @>>> \B U^q \,,
\ecd
where
$s_r$
and
$s_{r^q}$
denote the scalar multiplications by
$r$
and
$r^q \,.$

 Thus, the rational power of positive spaces emulates the rational
power of positive numbers, according to the above commutative
diagram.\END
\eNt

\bNt
 We have the following distinguished cases.

1) If
$q = 0 \,,$
then we have a natural semi--linear isomorphism
\beq
\B U^0 \byd \Rat^0(\B U^\star, \Rn^+) \seq \Rn^+ \,,
\eeq
and 
$\pi^0$ 
turns out to be the constant map with value 
$1 \,.$

\smallskip

2) Let
$q \eqv n$
be a positive integer.

 Then, we have the natural mutually inverse semi--linear isomorphisms
\bgt
\sten^n \B U \to \Rat^n(\B U^\star, \Rn^+) :
	u \sten \dots \sten u \to f_u \,,
\\
\Rat^n(\B U^\star, \Rn^+) \to \sten^n \B U :
f \mto u_f \sten \dots \sten u_f \,,
\end{gather*}
where
$f_u : \B U^\star \to \Rn^+ : \ome \mto \ome(u) \dots \ome(u)$
and
where
$u_f \byd 1/\ome_f \in \B U$
being 
$\ome_f \in \B U^\star$
the unique element such that
$f (\ome_f) = 1 \,.$

 Moreover, according to the above isomorphisms, the map
$\pi^n$ 
is given by
\beq
\pi^n : \B U \to \Rat^n(\B U^\star, \Rn^+) : u \mto f_u \,.
\eeq

 In particular, in the case
$q \eqv n = 1$
we have the natural semi--linear isomorphism
\beq
\B U^1 \byd \Rat^1(\B U^\star, \Rn^+)
= \Slin (\B U^\star, \Rn^+)
\byd \B U^{\star\star} \seq \B U \,.
\eeq

\smallskip

3) Let
$q \eqv 1/n$
be the inverse of a positive integer
$n \,.$

 Then, we have the natural mutually inverse semi--linear isomorphisms
\bgt
\sten^n \Rat^{1/n}(\B U^\star, \Rn^+) \to
\Slin (\B U^\star, \Rn^+)
\byd \B U^{\star\star} \seq \B U :
f \sten \dots \sten f \mto f^n \,,
\\
\B U \seq \B U^{\star\star} \byd \Slin (\B U^\star, \Rn^+) \to
\sten^n\Rat^{1/n}(\B U^\star, \Rn^+) :
f \mto f^{1/n} \sten \dots \sten f^{1/n} \,,
\end{gather*}
where
$f^n : \B U^\star \to \Rn^+ : \ome \mto f(\ome) \dots f(\ome)$
and
$f^{1/n} : \B U^\star \to \Rn^+ : \ome \mto (f(\ome))^{1/n} \,.$

 Moreover, according to the above isomorphisms, the map
$\pi^{1/n}$
is given by
\beq
\pi^{1/n} : \B U \seq \B U^{\star\star} \to 
\B U^{1/n} \byd \Rat^{1/n}(\B U^\star, \Rn^+) : f \mto f^{1/n} \,.
\eeq

\smallskip

4) Let
$q \eqv - n$
be a negative integer.
 
 Then, we have the natural mutually inverse semi--linear isomorphisms
\bgt
\sten^n \B U^\star \to \Rat^{-n}(\B U^\star, \Rn^+)
: \ome \sten \dots \sten \ome \to f_\ome \,,
\\
\Rat^{-n}(\B U^\star, \Rn^+) \to \sten^n \B U^\star :
f \mto \ome_f \sten \dots \ome_f \,,
\end{gather*}
where
$f_\ome : \B U^\star \to \Rn^+ : 
\alp \mto \ome (1/\alp) \dots \ome (1/\alp)$
and
$\ome_f \in \B U^\star$
is the unique element such that
$f(\ome_f) = 1 \,.$

 Moreover, according to the above isomorphisms, the map
$\pi^{-n}$ 
is given by
\beq
\pi^{-n} : \B U \to \Rat^{-n}(\B U^\star, \Rn^+) : u \mto f_{1/u} \,.
\eeq

 In particular, in the case
$q = -1 \,,$
we have the natural semi--linear isomorphism
\beq
\B U^{-1} \byd \Rat^{-1} (\B U^\star, \Rn^+)
\seq
\Slin(\B U^{\star\star}, \Rn^+) 
\seq
\Slin(\B U, \Rn^+) \byd
\B U^\star \,.\END
\eeq
\eNt

 Next, we analyse the natural behaviour of the exponents of rational
powers.
 Indeed, this behaviour is just what we expect and is analogous to
that of powers of positive real numbers.
 We leave to the reader the easy proofs of the following Propositions.

\bPr
 Let 
$p$ 
and 
$q$ 
be rational numbers.
 Then, we obtain the natural semi--bilinear map
\beq
b : \Rat^p (\B U^*, \Rn^+) \car \Rat^q (\B U^*, \Rn^+)
\to \Rat^{p+q}(\B U^*, \Rn^+) : (f, g) \mto fg \,,
\eeq
which yields the unique semi--linear isomorphism
$\ti b : \B U^p \sten \B U^q \to \B U^{p+q} \,,$
such that
$\pi^{p+q} = \ti b \com (\pi^p \sten \pi^q) \,,$
in virtue of the universal property of the semi--tensor product.\END
\ePr

\bPr
 If
$p$
and
$q$
are rational numbers, then we have the natural semi--linear
isomorphism
\beq
c : (\B U^p)^q \byd
\Rat^q \big(\Rat^p(\B U^*, \Rn^+), \Rn^+\big)
\to \B U^{p+q} \byd \Rat^{p+q}(\B U^*, \Rn^+) : f \mto g_f
\,,
\eeq
where
$g_f : \B U^\star \to \Rn^+ : 1/u \mto f(1/u^p) \,.$
 Moreover, we have
$\pi^{pq} = c \com (\pi^q \com \pi^p) \,.$\END
\ePr

\bCr
 If
$q$
is a rational number, then
\beq
(\B U^q)^* \seq (\B U^*)^q \,.
\eeq

 If 
$p < q$
are two positive integers, then
\bal
&\sten^q \B U \sten (\sten^p\B U^*) \seq
\B U^q \sten \B U^{-p} =
\B U^{q-p} \seq
\sten^{q-p}\B U
\\
&\sten^p \B U \sten (\sten^q\B U^*) \seq
\B U^p \sten \B U^{-q} =
\B U^{p-q} \seq
\sten^{|p-q|}\B U^* \,.\END
\end{align*}
\eCr
\section{Algebraic model of physical scales}
\label{Algebraic model of physical scales}
\bsm
 Next, we discuss the physical model of scales and units of
measurement through positive spaces.
 Thus, we introduce the fundamental scale spaces and related
notions, including scale dimension, scale basis and coupling scales.
 Finally, we show how unit spaces can be employed in physical theories
by the language of differential geometry.

 The formalism discussed in this section has been widely used in
several papers dealing with physical theories (see, for instance,
\cite{CanJadMod95,JanMod02c,JanMod05p2,ModSalTol05,SalVit00,Vit99,
Vit00}).
 In the present paper we analyse the mathematical foundations of this
formalism for the first time. We also discuss the interplay of our
theory with dimensional analysis.
\esm
\subsection{Units and scales}
\label{Units and scales}
\bsm
 We introduce the fundamental scale spaces and related notions.

 In this paper, we shall be concerned just with scales derived
from time, length and mass scales via rational powers.
 Of course, the treatment could be extended to other types of systems
in an analogous way.
\esm

 We assume the following positive spaces as \emp{basic spaces of
scales}:

(1) the space 
$\B T$ 
of \emp{time scales},

(2) the space 
$\B L$ 
of \emp{length scales},

(3) the space 
$\B M$ 
of \emp{mass scales}.

 The elements of the above spaces are called \emp{basic scales}.
 More precisely,

(1) each element
$u_0 \in \B T$
is said to be a \emp{time scale},

(2) each element
$\E l \in \B L$
is said to be a \emp{length scale},

(3) each element
$\E m \in \B M$
is said to be a \emp{mass scale}.

 For each time scale
$u_0 \in \B T \,,$
we denote its dual by
$u^0 \byd 1/u_0 \in \B T^* \,.$

\bDf
 A \emp{scale space} is defined to be a positive space of the type
\beq
\B S
\eqv \B S[d_1, d_2, d_3]
\byd \B T^{d_1} \sten \B L^{d_2} \sten \B M^{d_2} \,,
\ssep{where}
d_i \in \B Q \,.
\eeq

 A \emp{scale} is defined to be an element
$k \in \B S \,.$

 A scale
$k \in \B S \,,$
regarded as a semi--basis of the scale space
$\B S \,,$
is called a \emp{unit of measurement}.

 For each scale space
$\B S = \B T^{d_1} \sten \B L^{d_2} \sten \B M^{d_2}$
and for each scale
$k \in \B S \,,$
we set
\bal
\sdi {\B S}
&\eqv (\sdi {\B S}_1, \, \sdi {\B S}_2, \,\sdi {\B S}_3)
\byd (d_1, d_2, d_3)
\\
\sdi k
&\eqv (\sdi k_1, \, \sdi k_2, \,\sdi k_3)
\byd (d_1, d_2, d_3) \,.
\end{align*}

 The above 3-plet
$(d_1, d_2, d_3)$
of rational numbers is called the \emp{scale dimension} of
$\B S$
and of
$k \,.$\END
\eDf

 For instance, we have
$\B S[d_1, d_2, 0] \byd \B T^{d_1} \sten \B L^{d_2} \,,\,$
$\B S[d_1, 0, 0] \byd \B T^{d_1} \,,\,$
$\B S[0, 0, 0] \byd \Rn^+ \,,\,$
and the other similar examples.

 We stress that the scale dimension should not be confused with the
semi--dimension: indeed all these semi--vector spaces have
semi--dimension
$1 \,.$

\bNt
 The scale dimension
$\sdi k$
of a scale
$k$
determines the corresponding scale space
$\B S \,.$
 In other words, for two scales
$k$
and
$k' \,,$
we have
$\sdi k = \sdi{k'}$
if and only if the two scales belong to the same scale space. If this
is the case, then $k = r k' \,,$ where $r = k / k' \,.$
 Hence, the scale dimension
$\sdi k$
of a scale
$k$
determines
$k$
up to a positive real factor.

 The map
$k \mto \sdi k$
fulfills the following properties, for each
$k \in \B S \,,\,$
$k' \in \B S' \,,\,$
$r \in \Rn^+$
and
$q \in \B Q \,,$
\beq
\sdi {rk} = \sdi k \,,
\qquad
\sdi {1/k} = - \sdi k \,,
\qquad
\sdi {k \sten k'} = \sdi k + \sdi k' \,,
\qquad
\sdi {k^q} = q \, \sdi k \,.\END
\eeq
\eNt

\bDf
 A 3--plet of scales
$(e_1, e_2, e_3)$
is said to be a \emp{scale basis} if each scale
$k$
can be written in a unique way as
\beq
k = r \, (e_1)^{c_1} \sten (e_2)^{c_2} \sten (e_3)^{c_3} \,,
\ssep{with}
r \in \Rn^+ \,,
\;
c_i \in \B Q \,.\END
\eeq
\eDf

\bPr
 A 3--plet of scales
$(e_1, e_2, e_3)$
is a scale basis if and only if
\beq
\det (\sdi{e_j}_i) \neq 0 \,.
\eeq

 Moreover, let
$(e_1, e_2, e_3)$
be a scale basis and
$k$
a scale.
 Then, the 3--plet of rational exponents
$(c_1, c_2, c_3)$
is the unique solution of the linear rational system
\beq
\sdi{k}_i = \sum_j \sdi{e_j}_i \, c_j \,.
\eeq
\ePr

\bpf
 Let us consider a 3--plet of scales
$(e_1, e_2, e_3)$
and a scale
$k \,.$
 Then,
\beq
k = r \, (e_1)^{c_1} \sten (e_2)^{c_2} \sten (e_3)^{c_3}
\LRarr
\sdi{k}_i = \sum_j \sdi{e_j}_i \, c_j \,.
\eeq

 Hence, the above left hand side expression holds and is unique if and
only if
$\det (\sdi{e_j}_i) \neq 0 \,.$\QED
\epf

\bEx
 Clearly, each 3--plet of the type
\beq
(u_0, \E l, \E m) \in \B T \sten \B L \sten \B M
\eeq
is a scale basis.
 More generally, each 3--plet of the type
\beq
(u_0^{d_1}, \E l^{d_2}, \E m^{d_3}) \in
\B T^{d_1} \sten \B L^{d_2} \sten \B M^{d_3} \,,
\ssep{with}
d_i \in \B Q-\{0\} \,,
\eeq
is a scale basis.\END
\eEx

 Of course, we can also consider variable scales.
 Indeed, given a manifold
$\f M \,,$
we define a \emp{scale} of
$\f M$
to be a map of the type
$k : \f M \to \B S\,.$
\subsection{Scaled objects}
\label{Scaled objects}
\bsm
 In geometric models of physical theories one is often concerned with
vector bundle valued maps which have physical dimensions.
 Our theory of positive spaces allows us to keep into account this
fact in a formal algebraic way.
 In fact, we consider maps with values in vector bundles tensorialised
with positive spaces.
 The positive factors can be treated as numerical constants, with
respect to differential operators.
\esm

 Let
$\B U$
be a positive space. We observe that 
$\B U$ 
has a natural structure of
$1$-dimensional manifold.
 Moreover, it is easy to prove that the tangent space 
$T\B U$ 
is naturally isomorphic to a cartesian product.
 More precisely
\beq
T\B U \simeq \B U \times \baB U \,.
\eeq

Now, let us consider
a scale space
$\B S \,,$
two vector bundles
$p : \f F \to \f B$
and
$q : \f G \to \f B$
and a manifold
$\f M \,.$

 We can easily define the sesqui--tensor product bundle
$(\B U \lten \f F) \to \f B \,.$
 We can regard this vector bundle as the sesqui--tensor product over
$\f B$ 
of the trivial semi--vector bundle
$\wti{\B U} \byd (\f B \car \B U) \to \f B$
and the vector bundle
$\f F \to \f B \,.$

\bDf
 The bundle
$(\B S \lten \f F) \to \f B$
and its sections
$s : \f B \to \B S \lten \f F$
are said to be \emp{scaled}.
 Moreover, the bundle
$\B S \lten \f F \,,$ 
its sections
$s \in \sec(\f B, \, \B S \lten \f F) \,,$  
and the linear differential operator
\beq
\phi : \sec(\f B, \, \f G) \to \sec(\f B, \; \B S \lten \f F)
\eeq
are said to have \emp{scale dimension}
\beq
\sdi{\B S \lten \f F} = \sdi{s} = \sdi{\phi} = \sdi{\B S} \,.\END
\eeq
\eDf

\bNt
 Let 
$\C D : \sec(\f B, \f F) \to \sec(\f B, \f G)$
be a linear differential operator.
 Then, we obtain the linear differential operator (defined by the same
symbol)
\beq
\C D : 
\sec(\f B, \, \B S \lten \f F) \to \sec(\f B, \, \B S \lten \f G) :
s \mto \C D s \byd u \lten \C D \la \alp, 1/u \ra\,,
\eeq
where
$u \in \B S$
and
$\la \alp, 1/u \ra \in \sec(\f B, \f F) \,.$
 Of course, this definition does not depend on the choice of 
$u \,.$\END
\eNt

 The above construction applies, for instance, to the cases when
$\C D$
is the exterior differential, a Lie derivative, a covariant
derivative, and so on.

\bEx
 If
$\alp \in \sec(\f M, \, \B S \lten \Lam^r T^*\f M)$
is a scaled form, then we obtain the ``\emp{scaled exterior
differential}''
\beq
d \alp \byd u \lten d \alp' \in 
\sec(\f M, \, \B S \lten \Lam^{r+1} T^*\f M) \,,
\eeq
where
$u \in \B S$ 
and 
$\alp'$
is the form
$\alp' \byd \la \alp, 1/u \ra
\in \sec(\f M, \, \Lam^r T^*\f M) \,.$\END
\eEx

\bEx
 If
$t \in \sec(\f M, \, \ten^r T\f M)$
is a form and
$X \in \sec(\f M, \, \B S \lten T\f M)$
a scaled vector field. 
 Then, we obtain the ``\emp{scaled Lie derivative}''
\beq
L_X \, t \byd u \lten L_{X'} \, t \in 
\sec(\f M, \, \B S \lten (\ten^r T\f M)) \,,
\eeq
where
$u \in \B S$
and
$X'$
is the vector field
$X' \byd \la X, 1/u \ra  \in \sec(\f M,\,  T\f M) \,.$\END
\eEx

\bEx
 If
$c$
is a linear connection of the vector bundle
$\f F \to \f B \,,$
$X \in \sec(\f B, \, T\f B)$
a vector field and
$s \in \sec(\f B, \, \B S \lten \f F)$
a section.
 Then, we obtain the ``\emp{scaled covariant derivative}''
\beq
\nab_X s \byd u \lten \nab_X s' \in 
\sec(\f B, \, \B S \lten V\f F) \,,
\eeq
where
$u \in \B U$
and
$s'$
is the section
$s' \byd \la s, 1/u \ra \in \sec(\f B, \f F) \,.$

 We can re--interpret the above result in the following way.

 The trivial linear connection of
$(\f B \car \baB U) \to \f B$
and the linear connection
$c$
of
$\f F \to \f B$
yield a linear connection
$c'$
of
$(\B U \lten \f F) \to \f B \,,$
which has the same symbols of
$c \,.$
 The scaled covariant derivative can be regarded as the covariant
derivative with respect to the above product connection.

 By abuse of language, we shall denote by
$c$
also the product connection
$c' \,.$\END
\eEx
\subsection{Distinguished scales}
\label{Distinguished scales}
\bsm
 In this section, we discuss the algebraic model of distinguished
scales occurring in physics.
\esm

 Let us consider a vector bundle
$\f F \to \f B \,.$
 Suppose that in a physical theory we meet two scaled sections
with different scale factors
\beq
s : \f M \to \B S \lten \f F
\ssep{and}
s' : \f M \to \B S' \lten \f F \,.
\eeq

 Then, we can compare the two scales and write
$s = k \lten s' \,,$
provided we avail of a scale factor
$k : \f B \to \B S \sten \B S'{}^* \,,$
whose scale dimension is
$\sdi k = \sdi s - \sdi {s'} \,.$
 We call such a factor a \emp{coupling scale} (or, according to the
traditional use, a \emp{coupling constant}).

 Some coupling scales, such as, for instance, the speed of the
light, the Planck constant, the gravitational constant and the
positron charge have a fixed value, without reference to
specific systems.
 For this reason, we shall call these coupling scales \emp{universal}.

 Other types of coupling scales, such as, for instance, masses and
charges, arise, case by case, and are associated with different
particles.

 For instance, we have the following \emp{universal coupling scales}:

1) the \emp{speed of the light} 
$c \in \B T^{-1} \sten \B L \,,$

2) the  \emp{Planck constant}
$\h \in \B T^{-1} \sten \B L^2 \sten \B M \,,$

3) the \emp{gravitational constant}
$\E g \in \B T^{- 2} \sten \B L^3 \sten \B M^{-1} \,,$

4) the \emp{positron charge}
$\E e \in \B L^{3/2} \sten \B M^{1/2} \,.$

 Besides the above universal coupling scales, there are the following
scales which depend on the choice of a particle:

1) a \emp{mass}
$m \in \B M \,,$

2) a \emp{charge}
$q \in \baB T^{-1} \sten \B L^{3/2} \sten \B M^{1/2} \,.$
 We stress that a charge is a scale tensorialised with real numbers;
hence, a charge might be positive, vanishing, or negative.

\bNt
 The following 3--plets are scale bases (for
$q \neq 0$):
\bal
&1) \quad
(e_1, e_2, e_3) \byd (m, q, \h) \,,
\\
&2) \quad
(e_1, e_2, e_3) \byd (m, \h, \E g) \,,
\\
&3) \quad
(e_1, e_2, e_3) \byd (q, \h, \E g) \,.
\end{align*}

 Conversely, the following 3--plet is not a scale basis (for
$q \neq 0$):
\beq
4) \quad
(e_1, e_2, e_3) \byd (m, q, \E g) \,.
\eeq

 In fact, we have the following values of determinants in the above
cases, respectively:
\beq
1) \; \det (\sdi{e_j}_i) = - 1/2 \,,
\qquad
2) \; \det (\sdi{e_j}_i) = 1 \,,
\qquad
3) \; \det (\sdi{e_j}_i) = 1 \,,
\qquad
4) \; \det (\sdi{e_j}_i) = 0 \,.
\eeq

 Note that 
$\sdi{\E g} = \sdi{q^2/m^2} \,.$\END
\eNt

 It may be algebraically correct, but not physically reasonable to
express certain scales by means of some of the above scale bases.
 For instance, it may not be physically reasonable to express the
gravitational coupling scale
$\E g$
through the scale basis
$(m, q, \h) \,,$
because
$\E g$
is a universal coupling scale, while
$m$
and
$q$
depend on the choice of a specific particle.
\subsection{Interplay with dimensional analysis}
\label{dimensional analysis}
The \emph{dimensional analysis} (here, we use \cite{Bar03} as a
reference) is the branch of mathematical physics which studies the
properties of physical models which depend on units of measurement.

Many of the foundational ideas of dimensional analysis become very natural
facts in our algebraic theory of the units of measurement. Below, we list some
of these facts.
\begin{itemize}
\item A \emph{class of systems of units} \cite[p. 14]{Bar03} is, in
  our language, the choice of basic spaces of scales.
\item The \emph{dimension of a physical quantity} \cite[p. 16]{Bar03} is what
  we call the scale dimension. The \emph{dimension function of a physical
    quantity} is the expression of the physical quantity with respect to a
  given scale basis. Such an expression is always a rational function (provided
  that the quantity depends only on the chosen basic spaces of scales).
\item It can be proved \cite[p. 17]{Bar03} that the dimension function is
  always a power-law monomial. This justifies our algebraic setting: we obtain
  a rigorous formulation of these powers via tensor products and semi--linear
  duality. On the other hand, polynomials or power series would require
  additional constructions which seem not to be justified in view of the above
  property of the dimension function. Indeed, when in physics formulae
  containing power series occur they always involve real numbers, i.e. unscaled
  quantities (usually called ``pure numbers''), obtained as ratio of two scales
  belonging to the same positive space.  Often one of the two scales plays the
  role of a variable and the other one is regarded as a fixed distinguished
  scale.
\item The \emph{independence of dimensions} for some quantities \cite[p.\
  20]{Bar03} is just the property of those quantities of being a scale basis.
\item Any function that defines some relationship between quantities is
  homogeneous; this is a proposition from \cite[p. 24]{Bar03}.  It is a natural
  consequence of our setting that functions between scale spaces are rational.
  This property leads to the $\Pi$-theorem of dimensional analysis
  \cite{Bar03}.
\end{itemize}

Summarising, after realising that physical quantities transform with a
power-law monomial, it is natural to implement scales in physical models as
rational tensor powers, and maps between them as rational maps.

\textbf{Acknowledgements.} The last author would like to thank F. Paparella and
G. Saccomandi for useful discussions. This research has been supported by the
Ministry of Education of the Czech Republic under the project MSM0021622409,
Grant agency of the Czech Republic under the project GA 201/05/0523,
Universitiy of Florence, University of Salento, PRIN 2005-2007 ``Leggi di
conservazione e termodinamica in meccanica dei continui e in teorie di campo'',
GNFM and GNSAGA of INdAM.

\bfz

\efz
\end{document}